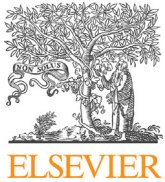
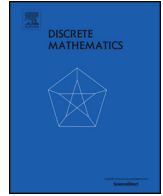

# New estimates for convex layer numbers

Gergely Ambrus [a,*], Peter Nielsen [b], Caledonia Wilson [c]

[a] *Alfréd Rényi Institute of Mathematics, Eötvös Loránd Research Network (ELKH), Reáltanoda u. 13-15, 1053 Budapest, Hungary*
[b] *Department of Mathematics, University of Wisconsin-Madison, 480 Lincoln Dr, Madison, WI 53706, USA*
[c] *Mount Holyoke College, 50 College St, MA 01075, USA*

### A B S T R A C T

Starting with a finite point set $X \subset \mathbb{R}^d$, the *peeling process* repeatedly removes the set of the vertices of the convex hull of the current set. The number of peeling steps required to completely remove $X$ is called the *layer number of $X$*, denoted by $L(X)$. In the article, we study the layer number of *evenly distributed families* of point sets contained in $B^d$, the $d$-dimensional unit ball. These sets consist of points in $B^d$ whose minimal distance is asymptotically as large as possible. We show that for a set $X$ belonging to an evenly distributed family, $L(X) \geq \Omega(|X|^{1/d})$ holds, with the bound being asymptotically sharp. On the other hand, building on earlier results, we prove that $L(X) \leq O(|X|^{2/d})$ holds for $d \geq 2$, which improves greatly on the current upper bound of $O(|X|^{(d+1)/2d})$ for $d \geq 3$. Finally, we provide a recursive construction of evenly distributed families whose sets satisfy $L(X) = \Theta(|X|^{2/d-1/(d2^{d-1})})$, showing that our upper bound is nearly tight.

© 2021 The Author(s). Published by Elsevier B.V. This is an open access article under the CC BY-NC-ND license (http://creativecommons.org/licenses/by-nc-nd/4.0/).

## 1. Introduction

Let $X$ be a finite point set in $\mathbb{R}^d$. The convex hull of $X$, denoted by conv$X$, is a convex polytope whose vertices are the *extreme points* of $X$. The set of extreme points will be denoted by $V(X)$. We call $V(X)$ the *first convex layer* of $X$.

Consider the following peeling process. Start with $X$. In each step, delete the extreme points of the current set from it – these constitute the *convex layers* of the set $X$. In finitely many steps, we reach the empty set. The number of steps is $L(X)$, the *layer number of $X$*. We note that in some sources, the layer number is referred to as the *convex depth* of $X$. Moreover, for $x \in X$, the *depth* of $x$ is the number of the step of the peeling process in which $x$ is deleted.

Layer numbers were first studied by Eddy [6] and Chazelle [4] from the algorithmic point of view. The latter article gave an optimal, $O(n \log n)$ running time algorithm for computing the convex layers of an $n$-element planar point set, showing that layer numbers may be computed quickly.

Almost 20 years later, Dalal [8] studied the layer number of random point sets. He proved that if $X$ is a set of $n$ random points chosen independently from the $d$-dimensional unit ball, then $L(X) = \Theta(n^{2/(d+1)})$. We note that studying the convex hull of random point sets (i.e. the first layer) dates back much earlier, see the works of Rényi and Sulanke [12] and Raynaud [11].

☆ Research of the first author was supported by NKFIH grants PD-125502 and KKP-133819. Work of the third author was supported by the Hungarian–American Fulbright Commission.
\* Corresponding author.
*E-mail addresses:* ambrus@renyi.hu (G. Ambrus), panielsen2@wisc.edu (P. Nielsen), wilso23c@mtholyoke.edu (C. Wilson).

https://doi.org/10.1016/j.disc.2021.112424

Har-Peled and Lidický [9] showed that the layer number of the planar $\sqrt{n} \times \sqrt{n}$ grid is $\Theta(n^{2/3})$. The analogous question for higher dimensions is studied in the recent article [2]. An interesting connection between the planar grid peeling process and the affine curve-shortening flow is presented in [7], while further algorithmic applications of the peeling process are given in [1].

In his Master's Thesis [10], W. Joo studied the layer numbers of $\alpha$-evenly distributed point sets in $\mathbb{R}^d$. For such sets, he proved the upper bound $L(X) \leq O(|X|^{\frac{d+1}{2d}})$. Moreover, he showed that this bound is sharp for $d = 2$. These results are published in the article [5] written jointly with I. Choi and M. Kim, which serves as the starting point for our current research.

In this paper, we are going to study evenly distributed families of point sets contained in $B^d$, the $d$-dimensional unit ball (see Definition 1.3 below). We prove that for sets $X$ belonging to such families in $\mathbb{R}^d$ with $d \geq 1$, $L(X) \geq \Omega(|X|^{\frac{1}{d}})$, and that this bound is sharp. On the other hand, we show that $L(X) \leq O(|X|^{\frac{2}{d}})$ (which improves on the current upper bound [5] of $O(|X|^{(d+1)/2d})$ for $d \geq 3$). For every $d \geq 2$ we construct point sets $X \subset B^d$ with $L(X) = \Theta(|X|^{(1-1/2^d)2/d})$ forming evenly distributed families, showing that our upper bound is nearly tight. Interestingly, the existence of an evenly distributed family of point sets $X$ in $\mathbb{R}^3$ satisfying $L(X) = \Theta(|X|^{\frac{2}{3}})$ would also lead to a construction with $L(X) = \Theta(|X|^{\frac{2}{d}})$ for all $d \geq 3$.

We start off with some basic notions and definitions. We are going to work in the $d$-dimensional Euclidean space $\mathbb{R}^d$ with the origin $o$. By distance, we refer to the Euclidean distance between two points, denoted by $|.|$. As usual, $B^d$ denotes the unit ball in $\mathbb{R}^d$, with $S^{d-1}$, the unit sphere in $\mathbb{R}^d$, being its boundary. As is well known,

$$s_{d-1} = d\kappa_d, \tag{1}$$

where $s_{d-1}$ denotes the surface area of $S^{d-1}$, and $\kappa_d$ gives the volume of $B^d$ (see [3]). For basic definitions regarding convexity we refer to [13].

The article [5] estimates layer numbers of $\alpha$-*evenly distributed* sets, which are defined as follows.

**Definition 1.1** ([5]). Let $X$ be a finite point set in a unit ball in $\mathbb{R}^d$. For a constant $\alpha > 1$, we say $X$ is $\alpha$-*evenly distributed* if for every Euclidean ball $D$ with positive volume,

$$|X \cap D| \leq \lceil \alpha |X| \text{Vol}(D) \rceil$$

For a set $A \subset \mathbb{R}^d$, let $\mu(A)$ denote the minimum (Euclidean) distance between points of $A$. The following assertion shows that being $\alpha$-evenly distributed is implied by a minimum distance condition.

**Lemma 1.2** (Lemma 2.3, [5]). *For every $d \geq 1$, there exists a continuous bijection $f_d : \mathbb{R}_{>0} \to \mathbb{R}_{>1}$ such that if $X \subset B^d$ is a finite point set satisfying $\mu(X) \geq \beta n^{-1/d}$ with some constant $\beta > 0$, then $X$ is $f_d(\beta)$-evenly distributed.*

In fact, being $\alpha$-evenly distributed is equivalent to the above property. In the current work, we study point sets satisfying this latter condition: their minimum distance is asymptotically as large as possible.

**Definition 1.3.** A family of sets $X_1, X_2, X_3, \cdots \subset B^d$ is said to be *evenly distributed* if $|X_i| \to \infty$ and $\mu(X_i) = \Theta(|X_i|^{-1/d})$.

The implied constant in $\Theta(|X_i|^{-1/d})$ may depend on $d$ as well as the family itself. Note that only the minimum distance between points of the sets in the family is guarded, yet, these sets may contain (large) holes.

By Lemma 1.2, every evenly distributed family is $\alpha$-evenly distributed for some parameter $\alpha$. In fact, the two definitions are equivalent to each other.

A key property that we use repeatedly is that if $X \subset B^d$ is a member of an evenly distributed family in $\mathbb{R}^d$ with $|X| = n$, then any ball of radius $O(n^{-1/d})$ contains at most $O(1)$ points of $X$. This follows by a standard volume estimate.

Most of the subsequent theorems regarding evenly distributed families have an asymptotic formulation. These are understood with the implied constants being dependent on the constant appearing in the asymptotic bound of Definition 1.3.

Our first results establish the sharp lower bound for evenly distributed families in $B^d$.

**Theorem 1.4.** *Assume that $\{X_i\}_{i=1}^{\infty} \subset B^d$ is an evenly distributed family in $\mathbb{R}^d$. Then $L(X_i) \geq \Omega(|X_i|^{1/d})$.*

**Proposition 1.5.** *For every $d \geq 1$, there exists an evenly distributed family $(X_i)_1^{\infty}$ with $L(X_i) = \Theta(|X_i|^{1/d})$.*

It is easy to see that for $d = 1$, the layer number of any set $X$ of distinct points in $\mathbb{R}$ is $\Theta(|X|)$. In the planar case, Choi, Joo, and Kim [5] proved that if $X$ is an $\alpha$-evenly distributed set with a parameter $\alpha > 1$, then $L(X) = O(|X|^{3/4})$, and this bound may be achieved. They also extended the bound to higher dimensions by showing that if $X$ is an $\alpha$-evenly distributed set in $\mathbb{R}^d$, then $L(X) = O(|X|^{(d+1)/2d})$, and they conjectured that this bound is in fact sharp. We improve their estimate to the following asymptotic upper bound:

**Theorem 1.6.** *Assume that $\{X_i\}_{i=1}^{\infty} \subset B^d$ is an evenly distributed family in $\mathbb{R}^d$. Then $L(X_i) \leq O(|X_i|^{2/d})$.*

Finally, we prove that the upper bound is asymptotically almost sharp.

**Theorem 1.7.** *For every $d \geq 1$, there exists an evenly distributed family $\{X_i^d\}_{i=1}^{\infty}$ in $\mathbb{R}^d$ with*

$$L(X_i^d) = \Theta\left(|X_i^d|^{\frac{2}{d} - \frac{1}{d2^{d-1}}}\right). \tag{2}$$

## 2. Well-separated sets

The quintessential example of evenly distributed families is given by the following standard definition.

**Definition 2.1.** A point set $X \subset \mathbb{R}^d$ is $\delta$-*separated* if $\mu(X) \geq \delta$. For $A \subset \mathbb{R}^d$ and $X \subset A$, we call $X$ to be *maximal $\delta$-separated in $A$* if $X$ is $\delta$-separated, and for any point $y \in A$, $|x - y| < \delta$ holds for some $x \in X$.

In other words, no further point of $A$ may be added to $X$ without losing the $\delta$-separated property. It is often used that maximal $\delta$-separated sets are also $\delta$-nets:

**Definition 2.2.** Let $A \subset \mathbb{R}^d$ and $X \subset A$. Then $X$ is a $\delta$-*net in $A$* if for every point $y \in A$, there exists $x \in X$ such that $|x - y| \leq \delta$.

We are going to make extensive use of the above notions with $A$ being the unit sphere $S^{d-1}$. In this case, the following well-known bound holds for the cardinality of maximal $\delta$-separated sets (see e.g. Section 2 of [3]).

**Lemma 2.3.** *Let $X$ be a maximal $\delta$-separated set in $S^{d-1}$. Then $X$ is a $\delta$-net, and*

$$|X| = \Theta(\delta^{-(d-1)}).$$

For self-containedness, we include the following sketch of the proof.

**Proof.** Since $X$ is a maximal $\delta$-separated set in $S^{d-1}$, for any point $y \in S^{d-1}$, $|x - y| < \delta$ holds for some $x \in X$. This shows that $X$ is indeed a $\delta$-net in $S^{d-1}$.

Additionally, since $X$ is $\delta$-separated set, Euclidean balls of radius $\delta/2$ centered at the points of $X$ are pairwise non-overlapping. For sufficiently small values of $\delta$, the intersection of such a ball with $S^{d-1}$ is well approximated by a spherical cap of radius $\delta/2$, which has surface area $\Theta(\delta^{d-1})$. Since these spherical caps are pairwise non-overlapping, we readily obtain that $|X| = O(\delta^{-(d-1)})$. On the other hand, since $X$ is a $\delta$-net in $S^{d-1}$, balls of radius $\delta$ centered at the points of $X$ cover $S^{d-1}$, and the same estimates show that $|X| = \Omega(\delta^{-(d-1)})$. $\square$

Clearly, the above asymptotic estimate remains true for any sphere with a fixed radius.

Without additional assumptions, the cardinality of a $\delta$-separated set in the unit ball $B^d$ may be as large as $\Theta(\delta^{-d})$. However, if the set is also in convex position, we may give a stronger estimate.

**Lemma 2.4.** *For any point set $X \subset B^d$ which is $\delta$-separated and in convex position,*

$$|X| \leq O(\delta^{-(d-1)}).$$

**Proof.** We show that the point set may be moved outwards to $S^{d-1}$ such that pairwise distances between points do not decrease, and then use the bound on the cardinality of a set in $S^{d-1}$ provided by Lemma 2.3.

Assume $x \in X$ is in the interior of $B^d$, that is, $|x| < 1$. Let $P = \text{conv} X$. Since $x$ is a boundary point of $P$, there exists an outer normal direction $u$ at $x$ (see e.g. [13]). That is, $u^{\perp} + x$ is a supporting hyperplane of $P$, with $u$ pointing away from $P$. Then for any positive $\lambda > 0$ and any $y \in X$, $y \neq x$, we have that $|y - x| < |y - (x + \lambda u)|$. Set $\lambda'$ so that $|x + \lambda' u| = 1$, and let $x' = x + \lambda' u$. Then if $X'$ is obtained from $X$ by replacing $x$ by $x'$, all pairwise distances of $X'$ are at least as large as the corresponding distance in $X$. By repeating the same process for each point of $X$ in the interior of $B^d$, we obtain a $\delta$-separated point set in $S^{d-1}$. The distance condition implies that distinct points remain distinct. By Lemma 2.3, the cardinality of this point set cannot exceed $\Theta(\delta^{-(d-1)})$. $\square$

## 3. The sharp lower bound on the layer number

In this section we establish our lower bound and prove its sharpness.

**Proof of Theorem 1.4.** By Definition 1.3, we know that $\mu(X_i) = \Theta(|X_i|^{-1/d})$. Thus, Lemma 2.4 implies that each convex layer of $X_i$ may have at most $O(|X|^{\frac{d-1}{d}})$ vertices. Hence, we obtain that $L(X) = \Omega(|X|^{1/d})$. □

It is not hard to prove that the above lower bound may be achieved.

**Proof of Proposition 1.5.** Let $S(1/2)$ denote the sphere of radius $1/2$ centered around the origin in $\mathbb{R}^d$. For $\lambda \geq 0$ and $A \subset \mathbb{R}^d$, we use the standard notation $\lambda A = \{\lambda a : a \in A\}$.

For every $i \geq 1$, let $\delta_i = 1/i$, and let $D_i$ be a maximal $\delta_i$-separated set on $S(1/2)$. Lemma 2.3 shows that $|D_i| = \Theta(\delta_i^{-(d-1)}) = \Theta(i^{(d-1)})$. Clearly, $D_i$ is in convex position.

To construct $X_i$, we let

$$X_i = \bigcup_{j=0}^{i-1} \left(1 + \frac{j}{i}\right) D_i.$$

Then, $|X_i| = i|D_i| = \Theta(i^d)$, and thus $\delta_i = \Theta(|X_i|^{-1/d})$. The shell structure of the construction implies that $L(X_i) = i = \Theta(|X_i|^{1/d})$. We only have to show that $\mu(X_i) = \Theta(|X_i|^{-1/d}) = \Theta(i^{-1})$ holds. The distance between points in the same layer of $X$ is, by definition, at least $\delta_i$. Finally, the distance between points on different layers is at least the difference between the radii of these layers, which is not less than $1/i = \delta_i$. Therefore, $X_i$ is $\delta_i$-separated, and the assertion follows. □

## 4. An upper bound on the layer number

The proof of Theorem 1.6 follows that in [5]. The key tool is the following statement therein.

**Lemma 4.1** ([5]). *If $K_1$ and $K_2$ are two convex bodies in $\mathbb{R}^d$ such that $K_2 \subseteq \text{int}(K_1)$ and $X$ is a finite point set in $K_1$, then*

$$L(X) \leq \max\{|X \cap C| : C \text{ is a cap of } K_1 \setminus K_2\} + L(X \cap K_2).$$

Here a *cap* of $K_1 \setminus K_2$ is the intersection $K_1 \cap H^+$, where $H$ is a supporting hyperplane of $K_2$ defining the corresponding closed halfspaces $H^-$ and $H^+$, of which $H^-$ contains $K_2$.

The proof of Lemma 4.1 may be found in [5]. The key idea is that if a cap $C$ satisfies $C \cap X \neq \emptyset$, then $C$ must also contain a point of $V(X)$. Thus, if $m$ denotes the maximal number of points of $X$ contained in a cap of $K_1 \setminus K_2$, then in $m$ steps of the peeling process, no point of $X$ may remain in any cap of $K_1 \setminus K_2$.

**Proof of Theorem 1.6.** Let $X$ be a member of an evenly distributed family in $B^d$ with $|X| = n$. Set $N = \lfloor n^{2/d} \rfloor$, and for $j \in \{0, \ldots, N\}$, let $B_j = (1 - \frac{j}{N})B^d$. For each $1 \leq j \leq N$, let $C_j$ be a cap of $B_{j-1} \setminus B_j$ which contains the maximum number of points of $X$ among all caps of $B_{j-1} \setminus B_j$. An elementary calculation reveals that the radius of the base of $C_j$ is at most

$$\sqrt{1 - (1 - N^{-1})^2} = O(n^{-1/d})$$

while its height is $\frac{1}{N} = o(n^{-1/d})$. Thus, $C_j$ is contained in a ball of radius $O(n^{-1/d})$. According to the remark following Definition 1.3, $|C_j \cap X| \leq O(1)$. Furthermore, since the radius of $B_N$ is at most $n^{-2/d} < n^{-1/d}$, we also have $|B_N \cap X| \leq O(1)$. Therefore, by Lemma 4.1,

$$L(X) \leq |X \cap B_N| + \sum_{j=0}^{N-1} |X \cap C_j|$$

$$\leq O(1) + \sum_{j=0}^{N-1} O(1)$$

$$\leq O(n^{\frac{2}{d}}). \quad \square$$

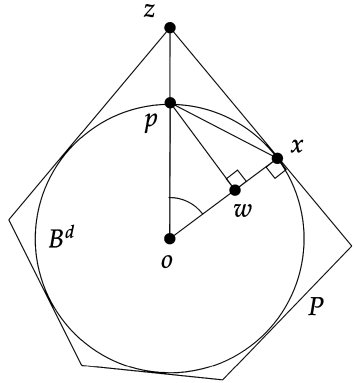

**Fig. 1.** We seek to find the smallest ball containing $P(X)$.

## 5. Tangent polytopes

The following notion plays a key role in the subsequent construction of point sets with large layer number.

**Definition 5.1.** Let $X \subset S^{d-1}$ be a finite point set affinely spanning $\mathbb{R}^d$ with $0 \in \mathrm{conv} X$. Define the *tangent polytope* $P(X)$ of $X$ as

$$P(X) = \cap_{x \in X} H^-(x)$$

where for each $x \in X$, $H(x)$ is the supporting hyperplane of $S^{d-1}$ at $x$, and $H^-(x)$ is the closed halfspace determined by $H(x)$ containing $S^{d-1}$.

The conditions $0 \in \mathrm{conv} X$ and $\dim(\mathrm{lin}\, X) = d$ ensure that $P(X)$ is a well-defined, bounded polytope in $\mathbb{R}^d$. For any $x \in X$, let $F(x)$ denote the (unique) face of $P(X)$ containing $x$.

Intuitively, the tangent polytope $P(X)$ may only be large if there is a big "gap" between points of $X$. The following statement formalizes this idea. Here and in the subsequent arguments, we assume that $\delta$ is small enough so that any $\delta$-net on $S^{d-1}$ contains the origin in its convex hull, and its points affinely span $\mathbb{R}^d$.

**Lemma 5.2.** *Let $X$ be a finite $\delta$-net in $S^{d-1}$ along with its tangent polytope $P = P(X)$. Then*

$$P \subset \frac{1}{1 - \delta^2/2} B^d.$$

**Proof.** Let $p \in S^{d-1}$ arbitrary, and let $z$ be the intersection point of the ray $op$ with the boundary of $P$. Then $z = \lambda p$ with some $\lambda \geq 1$. We are going to estimate $|z|$.

Let $x \in X$ be the point in $X$ closest to $p$. Since $X$ is a $\delta$-net in $S^{d-1}$, we have $|p - x| \leq \delta$.

As before, $H(x)$ is the supporting hyperplane of $S^{d-1}$ at $x$. Since the distance between $p$ and $H(x)$ is monotone increasing with respect to $|p - x|$, we necessarily have $z \in H(x)$. Consider the two-dimensional plane containing $z, x$, and $o$, see Fig. 1. Since $p \in H(x)$, $\angle oxz = \pi/2$.

Denote by $w$ the orthogonal projection of $p$ onto the segment $ox$, and let $|w| = t$. Applying the Pythagorean Theorem for $\triangle owp$ leads to $|p - w| = \sqrt{1 - t^2}$, and for $\triangle pwx$ yields

$$(1 - t)^2 + 1 - t^2 = |p - x|^2 \leq \delta^2.$$

Solving for $t$:

$$t \geq 1 - \frac{\delta^2}{2}.$$

Finally, since $\triangle owp \sim \triangle oxz$, we get that

$$|z| = \frac{|p|}{|w|}|x| = \frac{1}{t} \leq \frac{1}{1 - \delta^2/2}. \quad \square$$

We are also going to use the converse of Lemma 5.2.

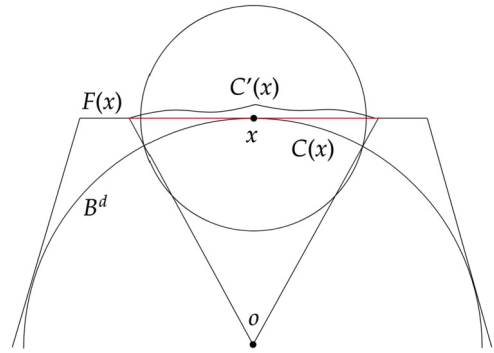

**Fig. 2.** Estimating the inradius of $F(x)$.

**Lemma 5.3.** *Let $X$ be a finite $\delta$-net in $S^{d-1}$, and set $P' = \text{conv} X$. Then*

$$\left(1 - \frac{\delta^2}{2}\right) B^d \subset P'. \tag{3}$$

**Proof.** Note that $P' = P(X)^\circ$, the polar body of the polytope $P(X)$, that is, $P^\circ = \{x \in \mathbb{R}^d : \langle x, y \rangle \leq 1 \text{ for every } y \in P(X)\}$ (see [13]). Since polarity reverses containment, and $(\lambda B^d)^\circ = (1/\lambda) B^d$, the statement follows. □

As usual, the *inradius* of a convex body in $\mathbb{R}^d$ is the radius of the largest ball contained within.

**Lemma 5.4.** *Let $X$ be a maximal $\delta$-separated set in $S^{d-1}$, with $P = P(X)$ being its tangent polytope. Then the inradius of each face of $P$ is at least $\delta/2$.*

**Proof.** By Lemma 2.3, $X$ is a $\delta$-net in $S^{d-1}$. Consider the spherical caps $C(x)$ of $S^{d-1}$ centered at points of $x \in X$ with (Euclidean) radius $\delta/2$:

$$C(x) = \left\{y \in S^{d-1} : |x - y| \leq \frac{\delta}{2}\right\}.$$

Since $X$ is a $\delta$-separated set, the triangle inequality shows that these spherical caps are pairwise disjoint.

For a fixed $x \in X$, let $H(x)$ be tangent hyperplane at $x$. Then $F(x) \subset H(x)$. Let $C'(x)$ be the radial projection of $C(x)$ onto $H(x)$ (that is, the set of intersection points between the hyperplane $H(x)$ and rays of the form $oy$ with $y \in C$). Since for every $y \in C(x)$, $x$ is the closest point to $y$ among points in $X$, we necessarily have that $C'(x) \subset F(x)$ (see Fig. 2). On the other hand, $C'(x)$ is a $(d-1)$-dimensional ball of radius greater than $\delta/2$. This completes the proof. □

## 6. Construction of evenly distributed sets with large layer numbers

Starting from $d = 1$, we are going to construct evenly distributed families of point sets in $B^d$ whose layers numbers are close to the upper bound provided by Theorem 1.6. The proof is motivated by the planar construction given in [5].

**Proof of Theorem 1.7.** We describe a recursive construction. For each $d \geq 1$, we are going to construct a family of sets $\{X_n^d\}_{n=1}^\infty$ which satisfy the following properties.

- *(P1)*: $X_n^d \subset B^d$
- *(P2)*: *For fixed $d$, $|X_n^d| = \Theta(n)$ with the implied constant depending on $d$*
- *(P3)*: *For fixed $d$, $\{X_n^d\}_{n=1}^\infty$ is evenly distributed in $\mathbb{R}^d$*
- *(P4)*: $o \in X_n^d$ *for every $d$ and $n$, and it is the last remaining point of the peeling process of $X_n^d$*
- *(P5)*: $L(X_n^d) = \Theta(n^{2/d - 1/(d 2^{d-1})})$.

We start with $d = 1$: for every $n \geq 1$, we set $D_n^1$ to be the set $\{i/n : i \in [-n, n]\}$. Then $D_n^1$ is a set of $2n + 1$ equally spaced points in $[-1, 1]$, whose layer number is $n + 1$. Therefore, it satisfies all the conditions *(P1) – (P5)*.

Now, let $d \geq 2$, and assume that the sets $X_n^{d-1}$ for every $n \geq 1$ have been constructed according to *(P1) – (P5)*. We are going to base our construction of $X_n^d$ on a parameter $\delta$, depending on $n$, whose value we will set later. Let us fix an arbitrary $n \geq 1$, and construct $X_n^d$.

Let $D$ be a maximal $\delta$-separated set in $S^{d-1}$ with tangent polytope $P = P(D)$. By Lemma 2.3,

$$|D| = \Theta(\delta^{-(d-1)}). \tag{4}$$

There are $|D|$ faces of $P$ of the form $F(x)$ with $x \in D$. Lemma 5.4 implies that within each face of $P$, a $(d-1)$-dimensional ball of radius $\delta/2$ may be inscribed with center $x$. Set

$$m = \lfloor \delta^{-(d-1)} \rfloor. \tag{5}$$

For every $x \in D$, embed in $F(x)$ a scaled copy of $X_m^{d-1}$ with scaling factor $\delta/4$ and center $x$. Denote the union of all these point sets on the faces of $P$ with $S$.

Let us study $S$. First, (4), (5) and property *(P2)* imply that

$$|S| = |D||X_m^{d-1}| = \Theta(\delta^{-2d+2}). \tag{6}$$

Second, we refer to the fact that $\mu(X_m^{d-1}) = \Theta(m^{-1/(d-1)})$, since $X_m^{d-1}$ is evenly distributed in $\mathbb{R}^{d-1}$ (see Definition 1.3). Thus, the minimum distance between points of $S$ contained in a given face $F(x)$ is

$$\frac{\delta}{4} \mu(X_m^{d-1}) = \Theta(\delta^2).$$

Because of disjointness of the $(d-1)$-dimensional balls of radius $\delta/2$ in the faces of $P$ centered at the points of $D$, and using the scaling factor $\delta/4$, the distance of a pair points of $S$ on different faces of $P$ is at least $\Omega(\delta) \gg \omega(\delta^2)$. Thus,

$$\mu(S) = \Theta(\delta^2). \tag{7}$$

Finally, since the sets $S \cap F(x)$ are congruent to each other (and $X_m^{d-1}$) for every $x \in D$, the $k$th layer removed by the peeling process is the union of the individual $k$th layers on the faces of $P$. By property *(P4)*, the last layer is the set $D$. Therefore, property *(P5)* implies that

$$L(S) = \Theta\left(m^{2/(d-1) - 1/((d-1)2^{d-2})}\right) = \Theta\left(\delta^{-2 + 1/2^{d-2}}\right). \tag{8}$$

Next, set

$$N = \left\lfloor \frac{1}{4\delta^2} \right\rfloor. \tag{9}$$

For every $i = 1, 2, \ldots, N$, let $r_i = 1 - i2\delta^2$, and let

$$S_i = r_i S. \tag{10}$$

Furthermore, let $\mathscr{S} = \cup_{i=1}^N S_i \cup \{o\}$. By (6),

$$|\mathscr{S}| = N|S| + 1 = \Theta(\delta^{-2d}). \tag{11}$$

The sets $S_i$ are to be called the *shells* of $\mathscr{S}$ (see Fig. 3). Note that $r_i \geq 1/2$ for every $i \leq N$. Thus, (7) implies that the minimum distance between points in the same shell of $S$ is $\Theta(\delta^2)$.

First, we show that $\mathscr{S} \subset B^d$. To that end, it is sufficient to prove that $(1 - 2\delta^2)S \subset B^d$. Since $S \subset P$, Lemma 5.2 shows that

$$S_1 = (1 - 2\delta^2)S \subset \frac{1 - 2\delta^2}{1 - \delta^2/2} B^d \subset (1 - \delta^2) B^d \tag{12}$$

which is clearly contained in $B^d$.

In order to establish the minimum distance condition of $\mathscr{S}$, we need to estimate the distance between different shells. Because of the scaling property, it is sufficient to give a lower bound on $|x - y|$ with $x \in S$ and $y \in S_1$: if $|x - y| \geq \theta$ holds for every such pair, then $|x' - y'| \geq \theta/2$ is true for every pair $x', y'$ contained in different shells. Clearly, $|x| \geq 1$ and, by (12), $|y| < 1 - \delta^2$. Thus, by the triangle inequality, $|x - y| > \delta^2$, and we obtain that

$$\mu(\mathscr{S}) = \Theta(\delta^2). \tag{13}$$

Next, we prove that

$$S_1 \subset \text{conv} D. \tag{14}$$

Indeed, by Lemma 5.3,

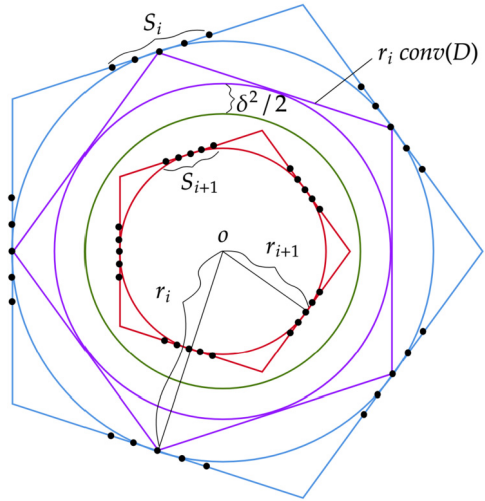

**Fig. 3.** Structure of the shell construction.

$$\left(1 - \frac{\delta^2}{2}\right) B^d \subset \text{conv} D.$$

Combined with (12), this shows that there is a ring of width $\delta^2/2$ separating $S_1$ and conv$D$.

By property *(P4)*, the points left by the penultimate step of the peeling process of each shell are the scaled copies of the set $D$. By applying (14), we obtain that the peeling process on $\mathscr{S}$ consists of the union of the individual peeling processes of the shells, one after the other. In addition, the last step of the peeling process of $\mathscr{S}$ is removing $\{o\}$. Thus, shells of $\mathscr{S}$ peel independently, and

$$L(\mathscr{S}) = N \cdot L(S) + 1 = \Theta\left(\delta^{-4+1/2^{d-2}}\right), \tag{15}$$

by (8) and (9).

Let us now set

$$\delta = n^{-1/(2d)}.$$

Then by (11), $\mathscr{S}$ is a point set with $\Theta(n)$ points, contained in $B^d$, which shows *(P1)* and *(P2)*. By (13), the minimum distance condition $\mu(\mathscr{S}) = \Theta(n^{-1/d})$ holds, implying *(P3)*. The origin $o$ is contained in $\mathscr{S}$, and it is the last point when peeling $\mathscr{S}$, which ensures that *(P4)* holds. Finally, by (15),

$$L(\mathscr{S}) = \Theta\left(n^{2/d - 1/d2^{d-1}}\right),$$

which agrees with *(P5)*. Thus, taking $X_n^d = \mathscr{S}$, the resulting family of sets satisfies all the conditions *(P1) – (P5)*. This concludes the proof. □

We remark that by defining $L_d$ so that $L(X_n^d) = \Theta(n^{L_d})$ in the above recursive construction, we obtain the following recurrence relation for $L_d$:

$$2dL_d = 2 + (d-1)L_{d-1} \tag{16}$$

with $L_1 = 1$. This may be solved by using exponential generating functions. Setting

$$F(x) = \sum_{k=1}^{\infty} x^k L_k,$$

(16) leads to

$$F'(x) = \frac{2}{2 - 3x + x^2}.$$

This, in turn, shows that

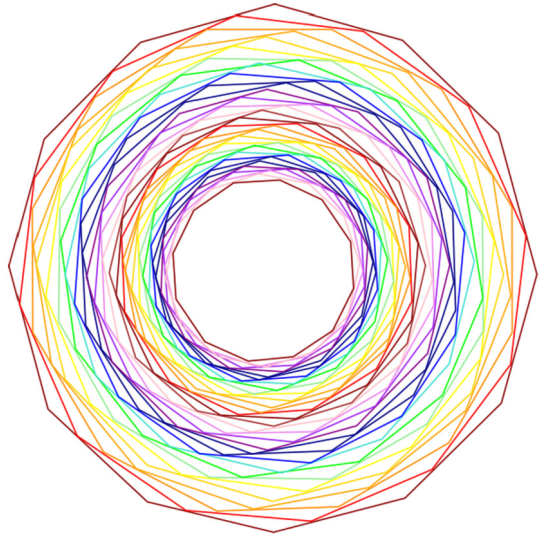

**Fig. 4.** Spiral construction in the plane.

$$L_d = \frac{2}{d} - \frac{2}{d2^d},$$

leading to the exponent in (2).

We also note that for $d = 2$, a different construction may be given with the same order of magnitude for the layer number. The construction consists of $n^{3/4}$ regular $n^{1/4}$-gons placed in a spiralling, interlocking manner (see Fig. 4). Given one layer, the next is obtained by moving $n^{-1/2}$ distance on each side of the polygon in positive orientation. Such a construction in $R^3$ gives a weaker bound than (2).

Finally we remark that the existence of an evenly distributed family of point sets $X$ in $\mathbb{R}^3$ satisfying $L(X) = \Theta(|X|^{\frac{2}{3}})$ would also lead to a construction with $L(X) = \Theta(|X|^{\frac{2}{d}})$ for all $d \geq 3$. Determining the exact asymptotics for the maximal layer number of an evenly distributed family remains an open question at this point.

### Declaration of competing interest

The authors declare that they have no known competing financial interests or personal relationships that could have appeared to influence the work reported in this paper.

### Acknowledgements

This research was done under the auspices of the Budapest Semesters in Mathematics program. We are grateful to W. Joo for communicating the results of [5] to us, and to the anonymous referees for their helpful suggestions.